\documentclass[11pt]{amsart}
\usepackage{amsmath,amsthm,amsfonts,amssymb,latexsym,epsfig}

\newtheorem{theorem}{Theorem}[section]

\newtheorem{lemma}{Lemma}[section]

\numberwithin{equation}{section}

\theoremstyle{definition}

\theoremstyle{remark}

\begin{document}
\title{Hardy-type Inequalities Via Auxiliary Sequences}
\author{Peng Gao}
\address{Department of Computer and Mathematical Sciences,
University of Toronto at Scarborough, 1265 Military Trail, Toronto
Ontario, Canada M1C 1A4} \email{penggao@utsc.utoronto.ca}
%%\date{\today}
\date{May 24, 2007.}
\subjclass[2000]{Primary 26D15} \keywords{Hardy's inequality}

%%-------------------------------------------------------------------------------------------

\begin{abstract}
  We prove some Hardy-type inequalities via an approach that
  involves constructing auxiliary sequences.
\end{abstract}

\maketitle
%%-------------------------------------------------------------------------------------------
\section{Introduction}
\label{sec 1} \setcounter{equation}{0}
%%-------------------------------------------------------------------------------------------

  Suppose throughout that $p\neq 0, \frac{1}{p}+\frac{1}{q}=1$.
   Let $l^p$ be the Banach space of all complex sequences ${\bf a}=(a_n)_{n \geq 1}$ with norm
\begin{equation*}
   ||{\bf a}||: =(\sum_{n=1}^{\infty}|a_n|^p)^{1/p} < \infty.
\end{equation*}
  The celebrated
   Hardy's inequality (\cite[Theorem 326]{HLP}) asserts that for $p>1$,
\begin{equation}
\label{eq:1} \sum^{\infty}_{n=1}\big{|}\frac {1}{n}
\sum^n_{k=1}a_k\big{|}^p \leq (\frac
{p}{p-1})^p\sum^\infty_{k=1}|a_k|^p.
\end{equation}
   Hardy's inequality can be regarded as a special case of the
   following inequality:
\begin{equation*}
\label{01}
   \sum^{\infty}_{j=1}\big{|}\sum^{\infty}_{k=1}c_{j,k}a_k
   \big{|}^p \leq U \sum^{\infty}_{k=1}|a_k|^p,
\end{equation*}
   in which $C=(c_{j,k})$ and the parameter $p$ are assumed
   fixed ($p>1$), and the estimate is to hold for all complex
   sequences ${\bf a}$. The $l^{p}$ operator norm of $C$ is
   then defined as the $p$-th root of the smallest value of the
   constant $U$:
\begin{equation*}
\label{02}
    ||C||_{p,p}=U^{\frac {1}{p}}.
\end{equation*}

    Hardy's inequality thus asserts that the Ces\'aro matrix
    operator $C$, given by $c_{j,k}=1/j , k\leq j$ and $0$
    otherwise, is bounded on {\it $l^p$} and has norm $\leq
    p/(p-1)$. (The norm is in fact $p/(p-1)$.)

    We say a matrix $A$ is a summability matrix if its entries satisfy:
    $a_{j,k} \geq 0$, $a_{j,k}=0$ for $k>j$ and
    $\sum^j_{k=1}a_{j,k}=1$. We say a summability matrix $A$ is a weighted
    mean matrix if its entries satisfy:
\begin{equation*}
%%\label{021}
    a_{j,k}=\lambda_k/\Lambda_j,  ~~ 1 \leq k \leq
    j; \Lambda_j=\sum^j_{i=1}\lambda_i, \lambda_i \geq 0, \lambda_1>0.
\end{equation*}

    Hardy's inequality \eqref{eq:1} now motivates one to
    determine the $l^{p}$ operator norm of an arbitrary summability matrix $A$.
   For examples, the following two
     inequalities were claimed to hold by Bennett ( \cite[p. 40-41]{B4}; see also \cite[p. 407]{B5}):
\begin{eqnarray}
\label{7}
   \sum^{\infty}_{n=1}\Big{|}\frac
1{n^{\alpha}}\sum^n_{i=1}(i^{\alpha}-(i-1)^{\alpha})a_i\Big{|}^p &
\leq & \Big( \frac {\alpha p}{\alpha p-1} \Big )^p\sum^{\infty}_{n=1}|a_n|^p, \\
\label{8}
   \sum^{\infty}_{n=1}\Big{|}\frac
1{\sum^n_{i=1}i^{\alpha-1}}\sum^n_{i=1}i^{\alpha-1}a_i\Big{|}^p &
\leq & \Big(\frac {\alpha p}{\alpha p-1} \Big
)^p\sum^{\infty}_{n=1}|a_n|^p,
\end{eqnarray}
     whenever $\alpha>0, p>1, \alpha p >1$.

     No proofs of the above two inequalities were supplied in \cite{B4}-\cite{B5} and
     recently, the author \cite{G} and Bennett himself \cite{Be1}
     proved inequalities \eqref{7} for $p>1, \alpha \geq 1, \alpha p >1$ and
     \eqref{8} for $p>1, \alpha \geq 2$ or $0< \alpha \leq 1, \alpha p >1$
     independently.

    We point out here that Bennett in fact was able
     to prove \eqref{7} for $p \geq 1, \alpha >0, \alpha p >1$ (see
     \cite[Theorem 1]{Be1} with $\beta=1$ there) which now leaves the
     case $p>1, 1< \alpha <2$ of inequality \eqref{8} the only case
     open to us. For this, Bennett expects inequality \eqref{8} to hold for $1+1/p < \alpha <2$
     (see page 830 of \cite{Be1}) and as a support,
     Bennett \cite[Theorem 18]{Be1} has shown that
inequality \eqref{8}
   holds for $\alpha = 1+1/p, p \geq 1$.

    In this paper, we will study inequality \eqref{8} using a
   method of Knopp \cite{K} which involves constructing auxiliary
   sequences. We will partially resolve the remaining case
   $p>1, 1< \alpha <2$ of inequality \eqref{8} by proving in Section
   \ref{sec 2} the following:
\begin{theorem}
\label{thm3}
   Inequality \eqref{8} holds for $p \geq 2, 1 \leq
   \alpha \leq 1+1/p$ or $1 < p \leq 4/3, 1+1/p \leq
   \alpha \leq 2$.
\end{theorem}

   We shall leave the explanation of Knopp's approach in detail in Section
   \ref{sec 2} by pointing out here that it can be applied to prove other types of inequalities
   similar to that of Hardy's. As an example, we note that Theorem
   359 of \cite{HLP} states:
\begin{theorem}
\label{thm4}
  For $0 < p <1$ and $a_n \geq 0$,
\begin{equation*}
  \sum^{\infty}_{n=1}\Big( \frac 1{n} \sum^{\infty}_{k=n}a_k \Big
  )^p \geq p^p \sum^{\infty}_{n=1}a^p_n.
\end{equation*}
\end{theorem}
  The constant $p^p$ in Theorem \ref{thm4} is not best possible
  and this was fixed by Levin and Ste\v ckin
  \cite[Theorem 61]{L&S} for $0<p \leq 1/3$ in the following
\begin{theorem}
\label{thm5}
  For $0 < p  \leq 1/3$ and $a_n \geq 0$,
\begin{equation*}
  \sum^{\infty}_{n=1}\Big( \frac 1{n} \sum^{\infty}_{k=n}a_k \Big
  )^p \geq \Big(\frac {p}{1-p} \Big
)^p \sum^{\infty}_{n=1}a^p_n.
\end{equation*}
\end{theorem}
   We shall give
  another proof of this result in Section \ref{sec 3}
  using Knopp's approach. We point out here for each $1/3< p <1$, Levin and Ste\v
  ckin also gave a better constant than the one $p^p$ given
  in Theorem \ref{thm4}. For example, when $p=1/2$, they
  gave $\sqrt{3}/2$ instead of $1/\sqrt{2}$. In Section \ref{sec 5}, we shall consider an approach of Redheffer \cite{R1} by showing first that this approach can be regarded as essentially the approach of Knopp when treating Hardy-type inequalities. 
We then use Redheffer's method to prove the following
\begin{theorem}
\label{thm6}
  For $a_n \geq 0$,
\begin{equation*}
  \sum^{\infty}_{n=1}\Big( \frac 1{n} \sum^{\infty}_{k=n}a_k \Big
  )^{1/2} \geq  0.8967 \sum^{\infty}_{n=1}a^{1/2}_n.
\end{equation*}
\end{theorem}
  This improves the result of Levin and Ste\v
  ckin mentioned above. It is also pointed out in Section \ref{sec 5} that the same method
can be used to establish the result in Theorem \ref{thm5} for $p$ slightly bigger than $1/3$.

  In our proofs of Theorems \ref{thm3}-\ref{thm4}, certain
  auxiliary sequences are constructed and there can be many ways to
  construct such sequences. In Section \ref{sec 4}, we give an
  example regarding these possibilities by answering a
  question of Bennett.

%%-------------------------------------------------------------------------------------
\section{Proof of Theorem \ref{thm3} }
\label{sec 2} \setcounter{equation}{0}
%%-------------------------------------------------------------------------------------
  We begin this section by explaining Knopp's idea \cite{K} on proving
  Hardy's inequality \eqref{eq:1}. In fact, we will explain this more generally for
  the case involving weighted mean matrices. For real numbers $\lambda_1>0, \lambda_i \geq 0,
  i \geq 2$, we write $\Lambda_n=\sum^n_{i=1}\lambda_i$ and we are
  looking for a positive constant $U$ such that
\begin{equation}
\label{2.0}
   \sum^{\infty}_{n=1}\Big{|} \frac {1}{\Lambda_n}\sum^{n}_{k=1}\lambda_ka_k
   \Big{|}^p \leq U \sum^{\infty}_{k=1}|a_k|^p
\end{equation}
   holds for all complex sequences ${\bf a}$ with $p>1$ being fixed.
   Knopp's idea is to find an
  auxiliary sequence ${\bf w}=\{ w_i \}^{\infty}_{i=1}$ of positive terms
  such that by H\"older's inequality,
\begin{eqnarray*}
\label{eq:12}
  \Big( \sum_{k=1}^n\lambda_k|a_k| \Big)^p &= & \Big( \sum_{k=1}^n\lambda_k|a_k|w_k^{-\frac {1}{p^*}} \cdot
  w_k^{\frac {1}{p^*}} \Big )^p \\
  & \leq & \Big( \sum_{k=1}^n\lambda^p_k|a_k|^pw_k^{-(p-1)} \Big ) \Big ( \sum_{j=1}^n w_j \Big)^{p-1}
\end{eqnarray*}
  so that
\begin{eqnarray*}\label{eq:13}
 \sum^{\infty}_{n=1}\Big{|} \frac {1}{\Lambda_n}\sum^{n}_{k=1}\lambda_ka_k
   \Big{|}^p  &\leq &
 \sum^\infty_{n=1}\frac
 {1}{\Lambda^p_n}\Big(\sum_{k=1}^{n}\lambda^p_k|a_k|^pw_k^{-(p-1)}\Big)\Big(\sum_{j=1}^n
 w_j\Big)^{p-1} \\
  &=& \sum_{k=1}^{\infty}w_k^{-(p-1)}\lambda^p_k\Big(\sum_{n=k}^{\infty}\frac
  {1}{\Lambda^p_n}\Big(\sum_{j=1}^nw_j\Big)^{p-1}\Big)|a_k|^p.
\end{eqnarray*}
  Suppose now one can find for each $p>1$ a positive constant $U$,
  a sequence ${\bf w}$ of positive terms with $w_n^{p-1}/\lambda^p_n$ decreasing to $0$, such
  that for any integer $n \geq 1$,
\begin{equation}
\label{eq:7}
 (w_1+\cdots+w_n)^{p-1}< U\Lambda^p_n( \frac {w_n^{p-1}}{\lambda^p_n}-\frac {w_{n+1}^{p-1}}{\lambda^p_{n+1}} ),
\end{equation}
  then it is easy to see that inequality \eqref{2.0} follows from this.
  When $\lambda_n=1$ for all $n$, Knopp's choice for ${\bf w}$ is given by
  $w_n=\binom{n-1-1/p}{n-1}$ and one can show that \eqref{eq:7}
  holds in this case with $U=(p^{*})^p$ and Hardy's inequality \eqref{eq:1}
  follows from this.

  We now want to apply Knopp's approach to prove Theorem \ref{thm3}. For this, we replace $\alpha-1$ by $\alpha$
  and rewrite \eqref{8} as
\begin{equation}
\label{2.10}
   \sum^{\infty}_{n=1}\Big{|}\frac
1{\sum^n_{i=1}i^{\alpha}}\sum^n_{i=1}i^{\alpha}a_i\Big{|}^p \leq
 \Big (\frac {(\alpha+1) p}{(\alpha +1)
 p-1} \Big )^p\sum^{\infty}_{n=1}|a_n|^p.
\end{equation}
  Note that we are interested in the case $0 \leq \alpha \leq 1$ here. From our
  discussions above, we are looking for a sequence ${\bf w}$ of positive terms
  with $w_n^{p-1}/\lambda^p_n$ decreasing to $0$, such
  that for any integer $n \geq 1$,
\begin{equation}
\label{2.20}
 (w_1+\cdots+w_n)^{p-1}< \Big (\frac {(\alpha+1) p}{(\alpha +1)
 p-1} \Big )^p\Big(\sum^n_{i=1}i^{\alpha} \Big)^p \Big ( \frac {w_n^{p-1}}{n^{\alpha p}}-\frac
{w_{n+1}^{p-1}}{(n+1)^{\alpha p}} \Big ).
\end{equation}
  Following Knopp's choice, we define a sequence ${\bf w}$ such that
\begin{equation}
\label{2.21}
  w_{n+1}= \frac {n+\alpha-1/p}{n}w_n, \hspace{0.1in} n \geq 1.
\end{equation}
   Note that the above sequence is uniquely determined for any given
   positive $w_1$ and therefore we may assume $w_1=1$ here. We note further that we need
   $\alpha > -1/p^{*}$ in order for $w_n >0$ for all $n$ and we also point
   out that it is easy to show by induction that
\begin{equation}
\label{2.22}
   \sum^{n}_{i=1}w_i=\frac {n+\alpha-1/p}{1+\alpha-1/p}w_n.
\end{equation}
  Moreover, one can easily check that
\begin{equation*}
  \frac {w_n^{p-1}}{n^{\alpha p}} = O(n^{-\alpha-1/p^{*}}),
\end{equation*}
   so that $w_n^{p-1}/\lambda^p_n$ decreases to $0$ as $n$
   approaches infinity as long as $\alpha > -1/p^{*}$.

    Now we need a lemma on sums of powers, which is due to Levin and Ste\v ckin
  \cite[Lemma 1, 2, p.18]{L&S}:
\begin{lemma}
\label{lem0}
    For an integer $n \geq 1$,
\begin{eqnarray}
\label{4}
    \sum^n_{i=1}i^r &\geq & \frac {1}{r+1}n(n+1)^r, \hspace{0.1in} 0 \leq r \leq 1, \\
\label{201}
   \sum^n_{i=1}i^r  &\geq & \frac {r}{r+1}\frac
   {n^r(n+1)^r}{(n+1)^r-n^r}, \hspace{0.1in} r \geq 1.
\end{eqnarray}
    Inequality \eqref{201} reverses when $-1 <r \leq 1$.
\end{lemma}
   We note here only the case $r \geq 0$ for \eqref{201} was proved in \cite{L&S} but one
   checks easily that the proof extends to the case $r >-1$.

   As we are interested in $0 \leq \alpha \leq 1$ here, we can now combine
  \eqref{2.21}-\eqref{4} to deduce that inequality \eqref{2.20} will
  follow from
\begin{equation*}
 (1+\frac {\alpha-1/p}{n})^{p-1} < \frac {n}{1+\alpha-1/p} \Big ((1+\frac 1{n})^{\alpha p}
 - (1+\frac {\alpha-1/p}{n} )^{p-1}\Big ).
\end{equation*}
  We can simplify the above inequality further by recasting it as
\begin{equation}
\label{2.23}
  \Big( 1+ \frac {\alpha+1/p^{*}}{n}\Big )^{1/p}\Big( 1+ \frac {\alpha-1/p}{n}\Big
  )^{1/p^{*}} < \Big( 1+ \frac {1}{n}\Big )^{\alpha}.
\end{equation}
  Now we define for
  fixed $n \geq 1, p>1$,
\begin{equation*}
  f(x)=x \ln (1+1/n)-\frac 1{p}\ln(1+ \frac {x+1/p^{*}}{n})-\frac 1{p^{*}}\ln(1+ \frac
  {x-1/p}{n}).
\end{equation*}
  It is easy to see here that inequality \eqref{2.23} is equivalent
  to $f(\alpha) > 0$. It is also easy to
  see that $f(x)$ is a convex function of $x$ for $0 \leq x \leq 1$
  and that $f(1/p)=0$. It follows from this that if $f'(1/p) \leq 0$
  then $f(x) > 0$ for $0 \leq x < 1/p$ and if $f'(1/p) \geq 0$ then
  $f(x) > 0$ for $1/p < x \leq 1$. We have
\begin{equation*}
  f'(1/p)=\ln (1+1/n)-\frac 1{n}+\frac 1{pn(n+1)}.
\end{equation*}
  We now use Taylor expansion to conclude for $x>0$,
\begin{equation}
\label{2.24}
  x - x^2/2< \ln (1+x) < x - x^2/2+ x^3/3.
\end{equation}
   It follows from this that for $p \geq 2$, $n \geq 2$,
\begin{equation*}
   f'(1/p)< -\frac 1{2n^2}+\frac 1{3n^3}+\frac 1{pn(n+1)} \leq
   -\frac 1{2n^2}+\frac 1{3n^3}+\frac 1{2n(n+1)}=\frac 1{3n^3}-\frac 1{2n^2(n+1)} \leq
   0.
\end{equation*}
   and for $n = 1$,
\begin{equation*}
   f'(1/p)= \ln 2-1 +\frac 1{2p} \leq  \ln 2-1 +\frac 1{4} < 0,
\end{equation*}

  It's also easy to check that for $1< p \leq 4/3$, $n=1$,
\begin{equation*}
  f'(1/p)=\ln 2-1+\frac 1{2p} > 0.
\end{equation*}
  For $n \geq 2, 1< p \leq 4/3$, by using the first inequality of
  \eqref{2.24} we get
\begin{equation*}
  f'(1/p) >  -\frac 1{2n^2}+\frac 1{pn(n+1)} \geq
   0.
\end{equation*}
  This now enables us to conclude the proof of Theorem \ref{thm3}.

%%-------------------------------------------------------------------------------------
\section{Another Proof of Theorem \ref{thm5}}
\label{sec 3} \setcounter{equation}{0}
%%-------------------------------------------------------------------------------------
   We use the idea of Levin and Ste\v ckin in the proof of Theorem
   62 in \cite{L&S} to find an auxiliary sequence ${\bf w}=\{ w_i
\}^{\infty}_{i=1}$ of positive
  terms so that for any finite summation from $n = 1$ to $N$ with $N \geq
  1$, we have
\begin{equation*}
  \sum^{N}_{n=1}a^p_n = \sum^{N}_{n=1}\frac {a^p_n}{\sum^n_{i=1}w_i} \sum^n_{k=1}w_k
  = \sum^N_{n=1}w_n\sum^{N}_{k=n}\frac {a^p_k}{\sum^k_{i=1}w_i}.
\end{equation*}
  On letting $N \rightarrow \infty$, we then have
\begin{equation*}
  \sum^{\infty}_{n=1}a^p_n =  \sum^{\infty}_{n=1}w_n\sum^{\infty}_{k=n}\frac {a^p_k}{\sum^k_{i=1}w_i}.
\end{equation*}
  By H\"older's inequality, we have
\begin{equation*}
  \sum^{\infty}_{k=n}\frac {a^p_k}{\sum^k_{i=1}w_i} \leq
  \Big ( \sum^{\infty}_{k=n}\Big(\sum^k_{i=1}w_i \Big)^{-1/(1-p)}\Big
  )^{1-p} \Big ( \sum^{\infty}_{k=n}a_k\Big
  )^{p}.
\end{equation*}
  Suppose now one can find
  a sequence ${\bf w}$ of positive terms with $w^{-1/(1-p)}_nn^{-p/(1-p)}$ decreasing to $0$ for each $0<p \leq 1/3$, such
  that for any integer $n \geq 1$,
\begin{equation}
\label{3.1}
 (w_1+\cdots+w_n)^{-1/(1-p)} \leq
 \Big ( \frac {1-p}{p} \Big )^{p/(1-p)}
 \Big( \frac {w^{-1/(1-p)}_n}{n^{p/(1-p)}}- \frac {w^{-1/(1-p)}_{n+1}}{(n+1)^{p/(1-p)}}\Big),
\end{equation}
   then it is easy to see that Theorem \ref{thm5} follows from this.

   We now define our sequence ${\bf w}$ to be
\begin{equation}
\label{3.2}
  w_{n+1}= \frac {n+1/p-2}{n}w_n, \hspace{0.1in} n \geq 1.
\end{equation}
   Note that the above sequence is uniquely determined for any given
   positive $w_1$ and therefore we may assume $w_1=1$ here. We note further that
   $w_n >0$ for all $n$ as $0 < p \leq 1/3$ and it is easy to show by induction that
\begin{equation}
\label{3.3}
   \sum^{n}_{i=1}w_i=\frac {n+1/p-2}{1/p-1}w_n.
\end{equation}
  Moreover, one can easily check that
\begin{equation*}
   \frac {w^{-1/(1-p)}_n}{n^{p/(1-p)}} = O(n^{-(1-p)/p}),
\end{equation*}
   so that $w_n^{-1/(p-1)}n^{-p/(1-p)}$ decreases to $0$ as $n$
   approaches infinity.

   We now combine
  \eqref{3.2}-\eqref{3.3} to recast inequality \eqref{3.1} as
\begin{equation*}
 (n+1/p-2)^{-1/(1-p)} \leq \frac {p}{1-p} \Big (n^{-p/(1-p)}
 - (n+1)^{-p/(1-p)}n^{1/(1-p)}(n+1/p-2)^{-1/(1-p)} \Big ).
\end{equation*}
   We further rewrite the above inequality as
\begin{eqnarray*}
   \frac {1-p}{p} & \leq & n^{-p/(1-p)}(n+1/p-2)^{1/(1-p)}
 - (n+1)^{-p/(1-p)}n^{1/(1-p)} \\
 &=& n \Big ( \Big(1+\frac {1/p-2}{n} \Big )^{1/(1-p)}
 - \Big (1+ \frac 1{n} \Big )^{-p/(1-p)} \Big ).
\end{eqnarray*}
  It is easy to see that the above inequality follows from $f(1/n) \geq 0$ where we define for $x \geq 0$,
\begin{equation*}
  f(x)=\Big(1+ (1/p-2)x \Big )^{1/(1-p)}
 - \Big (1+ x \Big )^{-p/(1-p)}-\frac {1-p}{p}x.
\end{equation*}
  We now prove that $f(x) \geq 0$ for $x \geq 0$ for $0 < p \leq 1/3$ and this will
  conclude the proof of Theorem \ref{thm5}. We note that
\begin{eqnarray*}
  f'(x) &=& \frac {1/p-2}{1-p}\Big(1+ (1/p-2)x \Big )^{p/(1-p)}+\frac p{1-p}\Big (1+ x \Big )^{-p/(1-p)-1}-\frac
  {1-p}{p},  \\
  f''(x) &=& \frac {p(1/p-2)^2}{(1-p)^2}\Big(1+ (1/p-2)x \Big )^{p/(1-p)-1}-\frac p{(1-p)^2}\Big (1+ x \Big
  )^{-p/(1-p)-2}.
\end{eqnarray*}
   We now define for $x \geq 0$,
\begin{equation*}
   g(x)=(1/p-2)^{2(1-p)/(1-2p)}(1+x)^{(2-p)/(1-2p)}-(1+ (1/p-2)x).
\end{equation*}
   It is easy to see that $g(x) \geq 0$ implies $f''(x) \geq 0$.
   Note that $(2-p)/(1-2p) \geq 1$ so that
\begin{eqnarray*}
   g'(x) &=&
   (1/p-2)^{2(1-p)/(1-2p)}(2-p)/(1-2p)(1+x)^{(2-p)/(1-2p)-1}-(1/p-2)
   \\
   & \geq & (1/p-2)^{2(1-p)/(1-2p)}-(1/p-2) \geq 0,
\end{eqnarray*}
   where the last inequality above follows from $2(1-p)/(1-2p) \geq
   1$ and $0 < p \leq 1/3$ so that $1/p - 2 \geq 1$. It follows from
   this that $f''(x) \geq 0$ and as one checks easily that
   $f'(0)=0$, which implies $f'(x) \geq 0$ so that $f(x) \geq f(0)=0$
   which is just what we want to prove.

%%-------------------------------------------------------------------------------------
\section{Redheffer's Approach and Proof of Theorem \ref{thm6} }
\label{sec 5} \setcounter{equation}{0}
%%-------------------------------------------------------------------------------------
  Redheffer's approach in  \cite{R1} of Hardy-type inequalities via his ``recurrent inequalities" can be put into the following form:  
\begin{lemma}[{\cite[Lemma 2.4]{G}}]
\label{lem6.1}
  Let $\{ \lambda_i \}^{\infty}_{i \geq 1}, \{ a_i \}^{\infty}_{i \geq 1}$ be two sequences of positive real numbers and
  let $S_n=\sum_{i=1}^n \lambda_i
  a_i$. Let $0 \neq p<1$ be fixed and let $\{ \mu_i \}^{\infty}_{i \geq 1}, \{ \eta_i \}^{\infty}_{i \geq 1}$ be two
  positive sequences of real numbers such
  that $\mu_i \leq \eta_i$ for $0<p<1$ and $\mu_i \geq \eta_i$ for
  $p<0$, then for $n \geq 2$,
\begin{equation}
\label{6.1}
   \sum_{i=2}^{n-1}\Big ( \mu_i-(\mu^q_{i+1}-\eta^q_{i+1})^{1/q} \Big )S_i^{1/p}+\mu_nS_n^{1/p}
   \leq (\mu^q_{2}-\eta^q_{2})^{1/q}\lambda^{1/p}_1a_1^{1/p}
   +\sum_{i=2}^n \eta_i \lambda^{1/p}_i a_i^{1/p}.
\end{equation}
\end{lemma}

   We consider the case $0<p<1$ in the above lemma and we set $\eta_i =
   \lambda^{-1/p}_i$ together with a change of variables: $\mu_i \mapsto \mu_i\eta_i$ to rewrite \eqref{6.1} as
\begin{equation*}
   \sum_{i=2}^{n-1}\Big ( \frac {\mu_i}{\lambda^{1/p}_i}-
   \frac {(\mu^q_{i+1}-1)^{1/q}}{\lambda^{1/p}_{i+1}} \Big )S_i^{1/p}+ \frac {\mu_n}{\lambda^{1/p}_n}S_n^{1/p}
   \leq (\mu^q_{2}-1)^{1/q}\frac {\lambda^{1/p}_1}{\lambda^{1/p}_2}a_1^{1/p}
   +\sum_{i=2}^n  a_i^{1/p}.
\end{equation*}
   We now set $\mu^q_{i}-1=\nu_i$ and make a further change of
   variables: $p \mapsto 1/p$ to write the above inequality as:
\begin{equation}
\label{6.2}
   \sum_{i=2}^{n-1}\Big ( \frac {(1+\nu_i)^{-(p-1)}}{\lambda^{p}_i}-
   \frac {\nu^{-(p-1)}_{i+1}}{\lambda^{p}_{i+1}} \Big )S_i^{p}+ \frac {(1+\nu_n)^{-(p-1)}}{\lambda^{p}_n}S_n^{p}
   \leq \nu^{-(p-1)}_2\frac {\lambda^{p}_1}{\lambda^{p}_2}a_1^{p}
   +\sum_{i=2}^n  a_i^{p}.
\end{equation}
   Now if we set for $i \geq 2$,
\begin{equation*}
   \nu_i=\frac {\sum^{i-1}_{j=1}w_j}{w_i},
\end{equation*}
   we can rewrite inequality \eqref{6.2} as
\begin{eqnarray}
\label{6.3}
  && \sum_{i=2}^{n-1}
    \Big(\sum^i_{j=1}w_j \Big )^{-(p-1)}\Big( \frac {w_i^{p-1}}{\lambda^p_i}-\frac {w_{i+1}^{p-1}}{\lambda^p_{i+1}} \Big )
    \Lambda^p_i A_i^{p}+
    \Big(\sum^n_{j=1}w_j \Big )^{-(p-1)}\frac {w_n^{p-1}}{\lambda^p_n}
    \Lambda^p_n A_n^{p}  \\
  &\leq & \frac {w_2^{p-1}}{w^{p-1}_1}\frac {\lambda^{p}_1}{\lambda^{p}_2}a_1^{p}
   +\sum_{i=2}^n  a_i^{p},  \nonumber
\end{eqnarray}
   where
\begin{equation*}
   \Lambda_n=\sum^n_{i=1}\lambda_i, \hspace{0.1in} A_n=\frac {S_n}{\Lambda_n}, \hspace{0.1in} n \geq 1.
\end{equation*}

   Suppose now we can find a sequence ${\bf w}=\{ w_i \}^{\infty}_{i=1}$ of positive
   terms such that inequality \eqref{eq:7} holds for all $n \geq 1$.
   Then inequality \eqref{6.3} implies
\begin{equation}
\label{6.4}
    \frac 1{U} \sum_{i=1}^{n} A_i^{p} \leq \Big ( \frac 1{U}+
    \frac {w_2^{p-1}}{w^{p-1}_1}\frac
    {\lambda^{p}_1}{\lambda^{p}_2} \Big ) a_1^{p}
   +\sum_{i=2}^n  a_i^{p} \leq \sum_{i=1}^n  a_i^{p},
\end{equation}
   where the last inequality above follows from the case $n=1$ of inequality \eqref{eq:7}, which implies
\begin{equation*}
   \frac 1{U} < 1-
    \frac {w_2^{p-1}}{w^{p-1}_1}\frac
    {\lambda^{p}_1}{\lambda^{p}_2}.
\end{equation*}
   Thus we have seen that on letting $n \rightarrow +\infty$,
   inequality \eqref{6.4} gives back inequality \eqref{2.0}.  Hence Redheffer's approach can be regarded as essentially Knopp's approach when treating Hardy-type inequalities. The only difference is that one no longer requires that $w_n^{p-1}/\lambda^p_n$ decreases to $0$ when selecting the sequence ${\bf w}$ in Redheffer's approach. 

%%-----------------------------------------------------------------------------------
%%-----------------------------------------------------------------------------------

%%------------------------------------------------------------------------------------
%%------------------------------------------------------------------------------------
   Now we state a lemma similar to Lemma \ref{lem6.1}:
\begin{lemma}
\label{lem6.2}
  Let $\{ \lambda_i \}^{\infty}_{i \geq 1}, \{ a_i \}^{\infty}_{i \geq 1}$ be two sequences of positive real numbers
  and suppose $\sum^{\infty}_{i=1}\lambda_ia_i$ converges.
  Let $S_n=\sum_{i=n}^{\infty} \lambda_i
  a_i$ and let $0 < p<1$ be fixed.  Let $\{ \mu_i \}^{\infty}_{i \geq 1}, \{ \eta_i \}^{\infty}_{i \geq 1}$ be two
  positive sequences of real numbers such
  that $\mu_i \geq \eta_i$, then for $n \geq 2$,
\begin{equation}
\label{6.5}
   \mu_1S^p_1+\sum_{i=2}^{n}\Big ( \mu_i-(\mu^{\frac 1{1-p}}_{i-1}-\eta^{\frac 1{1-p}}_{i-1})^{1-p} \Big )S_i^{p}
   -(\mu^{\frac 1{1-p}}_n-\eta^{\frac 1{1-p}}_n)^{1-p}S^{p}_{n+1}
   \geq \sum_{i=1}^n \eta_i \lambda^{p}_i a_i^{p}.
\end{equation}
\end{lemma}
\begin{proof}
   We note for
   $k \geq 2$,
\begin{equation}
\label{6.6}
    \mu_k S^{p}_k- \eta_k \lambda^{p}_k a_k^{p}=S^{p}_{k+1}(\mu_k (1+t)^{p}
    - \eta_kt^{p}) \geq (\mu^{\frac 1{1-p}}_k-\eta^{\frac 1{1-p}}_k)^{1-p}S^{p}_{k+1},
\end{equation}
     with $t=\lambda_k a_k/S_{k+1}$.
 The lemma then follows by adding \eqref{6.6} for $1 \leq k \leq n$
together.
\end{proof}

   We set $\eta_i =
   \lambda^{-p}_i$ together with a change of variables: $\mu_i \mapsto \mu_i\eta_i$ to rewrite \eqref{6.5} as
\begin{equation*}
   \frac {\mu_1}{\lambda^{p}_1}S_1^{p}+ \sum_{i=2}^{n}\Big ( \frac {\mu_i}{\lambda^{p}_i}-
   \frac {(\mu^{\frac 1{1-p}}_{i-1}-1)^{1-p}}{\lambda^{p}_{i-1}} \Big )S_i^{p}-
   \frac {(\mu^{\frac 1{1-p}}_{n}-1)^{1-p}}{\lambda^{p}_{n}} S_{n+1}^{p}
   \geq  \sum_{i=1}^n  a_i^{p}.
\end{equation*}
   We now set $\mu^{\frac 1{1-p}}_{i}-1=\nu_i$ and write the above inequality as:
\begin{equation*}
%%\label{6.7}
   \frac {(1+\nu_1)^{1-p}}{\lambda^{p}_1}S_1^{p}+\sum_{i=2}^{n}\Big ( \frac {(1+\nu_i)^{1-p}}{\lambda^{p}_i}-
   \frac {\nu^{1-p}_{i-1}}{\lambda^{p}_{i-1}} \Big )S_i^{p}- \frac {\nu_n^{1-p}}{\lambda^{p}_n}S_{n+1}^{p}
   \geq \sum_{i=1}^n  a_i^{p}.
\end{equation*}
   From now on we consider the case $\lambda_i=1$ for all $i$ in the above inequality and we set for $n \geq 1$,
\begin{equation*}
   \nu_n=\frac {n-\beta}{c},
\end{equation*}
   with $\beta \leq 1, c \geq \beta$ here. We want to choose $c, \beta$ such
   that the following inequality holds for $n \geq 2$:
\begin{equation}
\label{6.49}
  \max \Big( (1+c-\beta)^{1-p}, n^p\Big ((n+c-\beta)^{1-p}-(n-1-\beta)^{1-p} \Big ) \Big ) \leq  c^{1-p}k(p),
\end{equation}
  where $k(p)$ is a constant depending only on $p$ and we want
  $k(p)$ to be as small as possible.
  For this purpose, we further assume that $k(p)$ satisfies:
\begin{equation}
\label{6.50}
   (1-p)(1+c) < c^{1-p}k(p),
\end{equation}
   and define for $0 \leq x \leq 1/2$,
\begin{equation*}
    f(x) = (1+(c-\beta)x)^{1-p}-(1-(1+\beta)x)^{1-p}-c^{1-p}k(p)x,
\end{equation*}
   and note that with our assumption on $k(p)$, $f'(0)<0$. Note also
   that
\begin{equation*}
    f''(x)=p(1-p)(1+\beta)^2(1-(1+\beta)x)^{-p-1}-p(1-p)(c-\beta)^2(1+(c-\beta)x)^{-p-1}.
\end{equation*}
    It follows from this that when $1+\beta \geq c-\beta$ then
    $f''(x) \geq 0$ for $0 \leq x \leq 1/2$. Otherwise we note that
    $f''(x)=0$ can have at most one root in $(0, 1/2)$ and
    $f''(0)<0$. The above implies that for $0 \leq x \leq
    1/2$, $f(x) \leq \min (f(0), f(1/2))=\min (0, f(1/2))$.
    We deduce from our discussion above on setting $x=1/n$ in $f(x)$
    that in order for inequality \eqref{6.49} to hold, it suffices
    to check the case $n=2$, namely,
\begin{equation}
\label{6.51}
  \max \Big( (1+c-\beta)^{1-p}, 2^p\Big ((2+c-\beta)^{1-p}-(1-\beta)^{1-p} \Big ) \Big ) \leq  c^{1-p}k(p),
\end{equation}
   provided we assume \eqref{6.50}.

   We now look at the case $p=1/2$ and in this case we choose
   $c, \beta$ so that the following holds:
\begin{equation*}
  (1+\frac {1-\beta}{c})^{1/2}=2^{1/2}\Big ((1+\frac {2-\beta}{c})^{1/2}-(\frac {1-\beta}{c})^{1/2} \Big ).
\end{equation*}
   On setting $x=(1-\beta)/c, c'=1/c$, we can rewrite the above equation as
\begin{equation*}
   (1+x)^{1/2}=2^{1/2}\Big ( (1+c'+x)^{1/2}-x^{1/2} \Big ).
\end{equation*}
   Solving the above equation yields:
\begin{equation*}
   x=\frac {1-\beta}{c}=\frac {\sqrt{(10+4c')^2+28(1+2c')^2}-(10+4c')}{14}.
\end{equation*}
   To prove Theorem \ref{thm6}, we
   take $c=5/2$ here, then $x \approx 0.2435$ with $\beta \approx
   0.3912$ and $k(1/2)  \approx  1.1151 < 1.1152$. We take $k(p)=1.1152$ here and 
 one can also check that
   inequality \eqref{6.50} holds in this case. As $1/1.1152 > 0.8967$, Theorem \ref{thm6}
now follows.

   Now we consider other values of $p$'s. For this, we choose
   $c=1/p-1, k(p)=c^p$ so that inequality \eqref{6.50} becomes an
   equality. Because of this, we need to assume that
\begin{equation}
\label{6.54}
    \beta < \frac 1{2p}-1,
\end{equation}
   so that $f''(0)<0$ is satisfied for the function $f(x)$ defined above and our argument goes through as well in the above
   discussions to ensure that when inequality \eqref{6.51} holds,
   inequality \eqref{6.49} also holds. For the case $p=1/3$, it is easy to check that on taking $\beta = 3-2\sqrt{2}$,  both inequalities \eqref{6.51} and
   \eqref{6.54} are satisfied. This implies Theorem \ref{thm5} for the case $p =
   1/3$. In view of this, one sees that it is possible to prove the
   result in Theorem \ref{thm5} for $p$ beyond $1/3$. For example,
   on taking $p=0.34, \beta=0.21$, calculations shows both
   inequalities \eqref{6.54} and \eqref{6.51} are satisfied and
   hence Theorem \ref{thm5} holds for $p=0.34$.

%%-------------------------------------------------------------------------------------
\section{Another look at Inequality \eqref{8} }
\label{sec 4} \setcounter{equation}{0}
%%-------------------------------------------------------------------------------------
   In this section we return to the consideration of inequality \eqref{8} via our approach in
   Section \ref{sec 2}, which boils down to a construction of a sequence ${\bf w}$ of positive terms
  with $w_n^{p-1}/\lambda^p_n$ decreasing to $0$, such
  that for any integer $n \geq 1$, inequality \eqref{2.20} is
  satisfied. Certainly here the choice for ${\bf w}$ may not be
   unique and in fact in the case $\alpha =0$,
   Bennett asked in \cite{B4} (see the paragraph below Lemma 4.11) for other
  sequences, not multiples of Knopp's, that satisfy \eqref{2.20}. He
  also mentioned that the obvious choice, $w_n=n^{-1/p}$, does not
  work.

%%------------------------------------------------------------------------------
%%------------------------------------------------------------------------------
  We point out here even though the choice $w_n=n^{-1/p}$ does not
  satisfy \eqref{2.20} when $\alpha=0$ for all $p>1$,
  as one can see by considering inequality \eqref{2.20} for the case $n=1$
  with  $p \rightarrow 1^{+}$,
  it nevertheless works for $p \geq 3$, which we now show by first
  rewriting \eqref{2.20} in our case as
\begin{equation}
\label{2.30}
  \Big(\sum^n_{i=1}i^{-1/p} \Big )^{p-1}< \Big (\frac {p}{
 p-1} \Big )^pn^p\Big ( n^{-(p-1)/p }-
 (n+1)^{-(p-1)/p }
 \Big ).
\end{equation}
  We note that the case $n=1$ of \eqref{2.30} follows from the case $\alpha=0$ of the following inequality,
\begin{equation}
\label{2.4}
  1-2^{-(p-1)/p-\alpha} > \Big ( 1- \frac {1}{(\alpha+1)p} \Big )^p,
  \hspace{0.1in} 0 \leq \alpha \leq 1/p.
\end{equation}
  To show \eqref{2.4}, we see by Taylor expansion,
  that for $p \geq 2, x<0$,
\begin{equation*}
  (1+x)^p < 1+px+\frac {p(p-1)x^2}{2}.
\end{equation*}
  Apply the above inequality with $x=-1/(\alpha p+p)$, we obtain for $p \geq 3$,
\begin{equation*}
  \Big ( 1- \frac {1}{(\alpha+1)p} \Big )^p < 1-\frac {1}{(\alpha+1)}+\frac {(p-1)}{2(\alpha+1)^2p}.
\end{equation*}
  Hence inequality \eqref{2.4} will follow from
\begin{equation*}
  1-\frac {p-1}{2(\alpha+1)p}-2^{-(p-1)/p}\frac {(\alpha+1)}{2^{\alpha}} >0.
\end{equation*}
   It is easy to see that when $p \geq 3$, the function $\alpha \mapsto
   (1+\alpha)2^{-\alpha}$ is an increasing function of $\alpha$ for
   $0 \leq \alpha \leq 1/p$. It follows from this that for
   $0 \leq \alpha \leq 1/p$,
\begin{equation*}
  1-\frac {p-1}{2(\alpha+1)p}-2^{-(p-1)/p}\frac {(\alpha+1)}{2^{\alpha}}
  >  1-\frac {p-1}{2p}-2^{-(p-1)/p}\frac {(1/p+1)}{2^{1/p}} = 0,
\end{equation*}
  and from which inequality \eqref{2.4} follows.

  Now, to show \eqref{2.30} holds for all $n \geq 2, p \geq 3$, we
  first note that for $p>1$,
\begin{equation*}
   \sum^{n}_{i=1}i^{-1/p}  < 1+ \int^n_{1}x^{-1/p}dx = \frac {p}{p-1}n^{1-1/p}-\frac
   {1}{p-1}.
\end{equation*}

  On the other hand, by Hadamard's inequality,
  which asserts that for a continuous convex function $f(x)$ on $[a, b]$,
\begin{equation*}
   f(\frac {a+b}2) \leq \frac {1}{b-a}\int^b_a f(x)dx \leq \frac
   {f(a)+f(b)}{2},
\end{equation*}
   we have for $p>1$,
\begin{equation*}
   n^{-(p-1)/p}-(n+1)^{-(p-1)/p}=\frac {p-1}{p}\int^{n+1}_{n}x^{-1-1/p^{*}}dx
   \geq \frac {p-1}{p}(n+1/2)^{-1-1/p^{*}}.
\end{equation*}
  Hence inequality \eqref{2.30} will follow from the following
  inequality for $n \geq 2$,
\begin{equation*}
  \frac {p}{p-1}n^{1-1/p}-\frac
   {1}{p-1} \leq p^{*}n^{1/p^{*}}\Big ( 1 + \frac {1}{2n} \Big
   )^{-(1+p^{*})/p}.
\end{equation*}
   It is easy to see that for $p>1$,
\begin{equation*}
   \Big ( 1 + \frac {1}{2n} \Big
   )^{-(1+p^{*})/p} \geq 1 - \frac {1+p^{*}}{p}\frac {1}{2n}.
\end{equation*}
   Hence it suffices to show
\begin{equation*}
   \frac {p}{p-1}n^{1-1/p}-\frac
   {1}{p-1} \leq p^{*}n^{1/p^{*}}\Big ( 1 - \frac {1+p^{*}}{p}\frac {1}{2n} \Big
   ),
\end{equation*}
   or equivalently,
\begin{equation*}
  \Big ( 1+ \frac {1}{2p-2} \Big )^p \leq n.
\end{equation*}
  It's easy to check that the right-hand expression above is a
  decreasing function of $p \geq 3$ and is equal to $5^3/4^3<2$ when $p=3$.
  Hence it follows that \eqref{2.30} holds for all $n \geq 2, p \geq
  3$.
%%----------------------------------------------------------------------------
%%-----------------------------------------------------------------------------

   We consider lastly inequality \eqref{2.20} for other values of $\alpha$
   and we take $w_n=n^{\alpha-1/p}$ for $n \geq 1$ so that
  we can rewrite \eqref{2.20} as
\begin{equation}
\label{2.3}
  \Big(\sum^n_{i=1}i^{\alpha-1/p} \Big )^{p-1}< \Big (\frac {(\alpha+1) p}{(\alpha +1)
 p-1} \Big )^p\Big(\sum^n_{i=1}i^{\alpha} \Big)^p \Big ( n^{-(p-1)/p-\alpha }-
 (n+1)^{-(p-1)/p-\alpha }
 \Big ).
\end{equation}

   We end our discussion here by considering the case $1 \leq \alpha \leq 1+1/p$
   and we apply Lemma \ref{lem0} to obtain
\begin{eqnarray*}
  \sum^n_{i=1}i^{\alpha-1/p} &\leq & \frac {\alpha-1/p}{\alpha-1/p+1}\frac
   {n^{\alpha-1/p}(n+1)^{\alpha-1/p}}{(n+1)^{\alpha-1/p}-n^{\alpha-1/p}}
   =\frac {1}{\alpha-1/p+1}\Big ( \int^{n+1}_n x^{-\alpha+1/p-1}dx \Big
   )^{-1}, \\
  \sum^n_{i=1}i^{\alpha} &\geq & \frac {\alpha}{\alpha+1}\frac
   {n^{\alpha}(n+1)^{\alpha}}{(n+1)^{\alpha}-n^{\alpha}}
   =\frac {1}{\alpha+1}\Big ( \int^{n+1}_n x^{-\alpha-1}dx \Big
   )^{-1}
\end{eqnarray*}
  We further write
\begin{equation*}
   n^{-(p-1)/p-\alpha }-
 (n+1)^{-(p-1)/p-\alpha
 }=(\alpha-1/p+1)\int^{n+1}_{n}x^{-\alpha+1/p-2}dx,
\end{equation*}
  so that inequality \eqref{2.3} will follow from
\begin{equation*}
  \int^{n+1}_n x^{-\alpha-1}dx < \Big ( \int^{n+1}_n x^{-\alpha+1/p-1}dx \Big
   )^{1-1/p} \Big ( \int^{n+1}_{n}x^{-\alpha+1/p-2}dx \Big )^{1/p}.
\end{equation*}
   One can easily see that the above inequality holds by H\"older's
   inequality and it follows that inequality \eqref{2.3} holds for
   $p>1, 1 \leq \alpha \leq 1+1/p$. This provides another proof of
   inequality \eqref{8} for $p>1, 1 \leq \alpha \leq 1+1/p$.

%%-----------------------------------------------------------------------------------------


\begin{thebibliography}{99}
%%-----------------------------------------------------------------------------------------
%%\bibitem{al}
%%G. D. Allen, Power majorization and majorization of sequences, {\em
%%Result. Math.}, {\bf 14} (1988), 211-222.
%%----------------------------------------------------------------------------------------
%%\bibitem{alz}
%%H. Alzer, On an inequality of H. Minc and L. Sathre, {\em J. Math.
%%Anal. Appl.}, {\bf 179} (1993), 396-402.
%%----------------------------------------------------------------------------------------
%%\bibitem{alz1}
%%H. Alzer, Refinement of an inequality of G. Bennett, {\em  Discrete
%%Math.}, {\bf 135} (1994), 39-46.
%%-----------------------------------------------------------------------------------------
%%\bibitem{B&B} E. F. Beckenbach and R. Bellman, {\em
%%Inequalities}, Springer-Verlag,Berlin-G\"ottingen-Heidelberg 1961.
%%-----------------------------------------------------------------------------------------
%%\bibitem{Be0}
%%G. Bennett, Majorization versus power majorization, {\em Anal.
%%Math.}, {\bf 12} (1986), 283-286.
%%-----------------------------------------------------------------------------------------
%%\bibitem{Be}
%%G. Bennett, Lower bounds for matrices. II., {\em Canad. J. Math.},
%%{\bf 44} (1992), 54-74.
%%----------------------------------------------------------------
\bibitem{B4} G. Bennett, Factorizing the classical inequalities,
{\em Mem. Amer. Math. Soc.}, {\bf 120} (1996), 1--130.
%%----------------------------------------------------------------
\bibitem{B5} G. Bennett, Inequalities complimentary to Hardy, {\em
Quart. J. Math. Oxford Ser. (2)}, {\bf 49} (1998), 395--432.
%%-----------------------------------------------------------------------------------------
\bibitem{Be1}
G. Bennett, Sums of powers and the meaning of $l\sp p$, {\em Houston
J. Math.}, {\bf 32} (2006), 801-831.
%%-----------------------------------------------------------------------------------------
%%\bibitem{B&J} G. Bennett and G. Jameson, Monotonic averages of convex
%%functions, {\em J. Math. Anal. Appl.}, {\bf 252} (2000), 410-430.
%%----------------------------------------------------------------
%%\bibitem{Car} J. M. Cartlidge, Weighted mean matrices as operators on
%%$l^p$, Ph.D. thesis, Indiana University.
%%-----------------------------------------------------------------------------------------
%%\bibitem{CGQ}
%%T. H. Chan, P. Gao and F. Qi, On a generalization of Martin's
%%inequality, {\em  Monatsh. Math.}, {\bf 138} (2003), 179-187.
%%-----------------------------------------------------------------------------------------
%%\bibitem{Cl}
%%A. Clausing, A problem concerning majorization, in {\em General
%%Inequalities 4} (W. Walter, Ed.), Birkh\"user, Basel, 1984.
%%-----------------------------------------------------------------------------------------
\bibitem{G}
P. Gao, {A note on Hardy-type inequalities}, {\it Proc. Amer. Math.
Soc.}, \textbf{133} (2005), 1977-1984.
%%-----------------------------------------------------------------------------------------
%%\bibitem{G1}
%%P. Gao, {On a result of Cartlidge}, {\it J. Math. Anal. Appl.},
%%accepted.
%%-----------------------------------------------------------------------------------------
\bibitem{HLP} G. H. Hardy, J. E. Littlewood and G. P\'{o}lya, {\em
Inequalities}, Cambridge Univ. Press, 1952.
%%-----------------------------------------------------------------------------------------
\bibitem{K} K. Knopp, \"Uber Reihen mit positiven Gliedern,
\emph{J. London Math. Soc.}, {\bf 3} 1928, 205-211 and {\bf 5} 1930,
13-21.
%%----------------------------------------------------------------------------------------
\bibitem{L&S} V. I. Levin and S.B. Ste\v ckin, Inequalities, {\em
Amer. Math. Soc. Transl. (2)}, {\bf 14} (1960),  1--29.
%%-----------------------------------------------------------------------------------------
\bibitem{R1} R. M. Redheffer, Recurrent inequalities, {\em Proc.
London Math. Soc. (3)}, {\bf 17} (1967), 683--699.
%%-----------------------------------------------------------------------
%%\bibitem{sto} K. B. Stolarsky, Generalizations of the
%%logarithmic mean, {\em Math. Mag.}, {\bf 48} (1975), 87-92.
%%----------------------------------------------------------------
%%\bibitem{M&O&P}
%%A. W. Marshall, I. Olkin and F. Proschan, Monotonicity of ratios of
%%means and other applications of majorization, {\em Inequalities
%%(Proc. Sympos. Wright-Patterson Air Force Base, Ohio, 1965)}, pp.
%%177-190, Academic Press, New York, 1967.
%%-----------------------------------------------------------------------------------------
%%\bibitem{Mar}
%%J. S. Martins, Arithmetic and geometric means, an application to
%%Lorentz sequence spaces, {\em Math. Nachr.}, {\bf 139} (1988),
%%281--288.
%%-----------------------------------------------------------------------------------------
%%\bibitem{Xu}
%%Z. K. Xu and D. P. Xu, A general form of Alzer's inequality, {\em
%%Comput. Math. Appl.}, {\bf 44} (2002), 365-373.


\end{thebibliography}
\end{document}